\documentclass[a4paper,11pt]{amsart}
\usepackage{wrapfig}
\usepackage[dvips]{graphicx}
\usepackage{amsmath,amsthm,amssymb,amscd}
\theoremstyle{definition}
\newtheorem{thm}{Theorem}[section]
\newtheorem{Def}[thm]{Definition}
\newtheorem{pro}[thm]{Proposition}
\newtheorem{cor}[thm]{Corollary}
\newtheorem{lem}[thm]{Lemma}
\newtheorem{ex}[thm]{Example}
\newtheorem{rem}[thm]{Remark}
\theoremstyle{definition}

\begin{document}

\title{Fundamental group of simple $C^*$-algebras with unique trace}
\author{Norio Nawata}
\address[Norio Nawata]{Graduate School of Mathematics, 
Kyushu University, Hakozaki, 
Fukuoka, 812-8581,  Japan}      

\author{Yasuo Watatani}
\address[Yasuo Watatani]{Department of Mathematical Sciences, 
Kyushu University, Hakozaki, 
Fukuoka, 812-8581,  Japan}
\maketitle
\begin{abstract}
We introduce the fundamental group ${\mathcal F}(A)$ of 
a unital simple $C^*$-algebra $A$ with a unique normalized 
trace.  We compute fundamental groups 
${\mathcal F}(A)$ of 
several  nuclear or non-nuclear $C^*$-algebras $A$. 
K-theoretical obstruction enables us to compute the 
fundamental group easily. 
Our study is essentially based on the computation  of 
Picard groups by Kodaka.

\end{abstract}

\section{Introduction}
Let $M$ be a factor of type $II_1$ with a normalized trace 
$\tau$. Murray and von Neumann  introduced 
the fundamental group ${\mathcal F}(M)$ of $M$ in \cite{MN}. 
They showed that if $M$ is  hyperfinite, then 
${\mathcal F}(M) = {\mathbb R_+^{\times}}$. Since then 
there has been many works on the computation of the 
fundamental groups. Voiculescu \cite{Vo} showed that 
${\mathcal F}(L(F_{\infty}))$ of the group factor 
of the free group $F_{\infty}$ contains the positive rationals and 
Radulescu proved that 
${\mathcal F}(L(F_{\infty})) = {\mathbb R}_+^{\times}$ in 
\cite{Ra}.  Connes \cite{Co} showed that if $G$ is an ICC group with property 
(T), then  ${\mathcal F}(L(G))$ is a countable group. Popa 
showed that any countable subgroup of $\mathbb R_+^{\times}$ 
can be realized as the fundamental group of some 
factor of type $II_1$ in \cite{Po1}. 

In this paper  
we introduce the fundamental group of simple $C^*$-algebras 
with unique trace. 
Our study is essentially based on the computation  of 
Picard groups by Kodaka  \cite{kod1}, \cite{kod2}, \cite{kod3}.  
Let $A$ be a unital simple $C^*$-algebra with a unique normalized 
trace $\tau$. The fundamental group ${\mathcal F}(A)$ of $A$ is 
defined as the set of 
the numbers $\tau \otimes Tr(p)$ for some projection 
$p \in M_n(A)$ such that 
$pM_n(A)p$ is isomorphic to $A$. Then 
the fundamental group ${\mathcal F}(A)$ of $A$ is a 
multiplicative subgroup of ${{\mathbb R}_+^{\times}}$. 
After studying  basic properties of the fundamental groups of 
unital simple $C^*$-algebras, we compute ${\mathcal F}(A)$ of 
several nuclear or non-nuclear $C^*$-algebras $A$. 
If a $C^*$-algebra $A$ belongs to the class covered by Elliott's 
classification theorem, then the fundamental group ${\mathcal F}(A)$ 
can be computed rather easily.  But even if $A$ is non-nuclear, 
numerical invariant  ${\mathcal F}(A)$ is often computable. 
For example, 
let $A_{\theta}$ be an irrational rotation algebra, then 
${\mathcal F}(A_{\theta})$ is trivial, if $\theta$ is not 
a quadratic number, and ${\mathcal F}(A_{\theta})$ 
is the group of positive invertible elements  of 
the ring 
${\mathbb Z} + {\mathbb Z}\theta$ if 
$\theta$ is 
a quadratic algebraic integer. For example, 
${\mathcal F}(A_{\sqrt{3}}) 
= \{(2 + \sqrt{3})^n | n \in {\mathbb Z} \}.$
For a free group $F_n$ of $n$ generators with $n \geq 2$,  
we have ${\mathcal F}(A_{\theta} \otimes C_r^*(F_n)) 
= {\mathcal F}(A_{\theta})$. Some computations suggest that 
the fundamental groups 
of simple $C^*$-algebras have arithmetical flavour  
compared with the case of factors of type $II_1$.

\section{Hilbert $C^*$-modules and Picard groups}\label{sec:picard}

Let $A$ be a $C^*$-algebra and 
$\mathcal{E}$ a right Hilbert $A$-module. We denote by $L_A(\mathcal{E})$ 
the algebra of the adjointable operators on $\mathcal{E}$. 
For $\xi,\eta \in \mathcal{E}$, a  "rank one operator" $\Theta_{\xi,\eta}$
is defined by $\Theta_{\xi,\eta}(\zeta) 
= \xi \langle\eta,\zeta\rangle_A$ for $\zeta \in \mathcal{E}$. 
We denote by $K_A(\mathcal{E})$ the closure 
of the linear span of "rank one operators" $\Theta_{\xi,\eta}$.
We call a finite 
set $\{\xi_i\}_{i=1}^n\subseteq \mathcal{E}$ a {\it finite basis} of $\mathcal{E}$ if 
$\eta =\sum_{i=1}^n\xi_i\langle\xi_i,\eta\rangle_A$ for any $\eta\in\mathcal{E}$, see \cite{KW}, \cite{W}. It is also called a frame as in \cite{FL}.   
If $A$ is unital and there exists a finite basis for  $\mathcal{E}$, then $L_A(\mathcal{E})=K_A(\mathcal{E})$. 

We recall some definitions and elementary facts on the Picard groups of 
$C^*$-algebras introduced by Brown, Green and Rieffel 
in \cite{BGR}. 
An $A$-$A$-equivalence bimodule is an $A$-$A$-bimodule $\mathcal{E}$ which is simultaneously a 
full left Hilbert $A$-module under a left $A$-valued inner product $_A\langle\cdot ,\cdot\rangle$ 
and a full right Hilbert $A$-module under a right $A$-valued inner product $\langle\cdot ,\cdot\rangle_A$, 
satisfying $_A\langle\xi ,\eta\rangle\zeta =\xi\langle\eta ,\zeta\rangle_A$ for any 
$\xi, \eta, \zeta \in \mathcal{E}$. For $A$-$A$-equivalence bimodules 
$\mathcal{E}_1$ and 
$\mathcal{E}_2$, we say that $\mathcal{E}_1$ is isomorphic to $\mathcal{E}_2$ as an equivalence 
bimodule if there exists a $\mathbb{C}$-liner one-to-one map $\Phi$ from $\mathcal{E}_1$ onto 
$\mathcal{E}_2$ with the properties such that $\Phi (a\xi b)=a\Phi (\xi )b$, 
$_A\langle \Phi (\xi ) ,\Phi(\eta )\rangle =\;_A\langle \xi ,\eta\rangle$ and 
$\langle \Phi (\xi ) ,\Phi(\eta )\rangle_A =\langle\xi,\eta\rangle_A$ for $a,b\in A$, 
$\xi ,\eta\in\mathcal{E}_1$. 
The set of isomorphic classes $[\mathcal{E}]$ of the $A$-$A$-equivalence 
bimodules $\mathcal{E}$ forms a group under the product defined by 
$[\mathcal{E}_1][\mathcal{E}_2] = [\mathcal{E}_1 \otimes_A \mathcal{E}_2]$. 
We call it the {\it Picard group} of $A$ and denote it  by  $\mathrm{Pic}(A)$. 
The identity of $\mathrm{Pic}(A)$ is given by 
the $A$-$A$-bimodule $\mathcal{E}:= A$ with  
$\; _A\langle a ,b \rangle = ab^*$ and $\langle a ,b\rangle_A = a^*b$ for 
$a,b \in A$. 
The dual module $\mathcal{E}^*$ of an $A$-$A$-equivalence bimodule $\mathcal{E}$ is a set 
$\{\xi^* ;\xi\in\mathcal{E} \}$ with the operations such that $\xi^* +\eta^*=(\xi +\eta )^*$, 
$\lambda\xi ^*=(\overline{\lambda}\xi)^*$, $b\xi^* a=(a^*\xi b^*)^*$, 
$_A\langle\xi^*,\eta^*\rangle =\langle\eta ,\xi\rangle_A$ and 
$\langle \xi^*,\eta^*\rangle_A =\;_A\langle\eta ,\xi\rangle$. Then  
 $[\mathcal{E}^*]$ is the  
inverse element of $[\mathcal{E}]$ in the Picard group of $A$. 
Let $\alpha$ be an automorphism of $A$, and let 
$\mathcal{E}_{\alpha}=A$ with the obvious left $A$-action and the obvious $A$-valued inner product. 
We define the right $A$-action on $\mathcal{E}_\alpha$ by 
$\xi\cdot a=\xi\alpha(a)$ for 
any $\xi\in\mathcal{E}_\alpha$ and $a\in A$, and the right $A$-valued inner product by 
$\langle\xi ,\eta\rangle_A=\alpha^{-1} (\xi^*\eta)$ for any $\xi ,\eta\in\mathcal{E}_\alpha$.
Then $\mathcal{E}_{\alpha}$ is an $A$-$A$-equivalence bimodule. For $\alpha, \beta\in\mathrm{Aut}(A)$, 
$\mathcal{E}_\alpha$ is isomorphic to $\mathcal{E}_\beta$ if and only if 
there exists a unitary $u \in A$ such that 
$\alpha = ad \ u \circ \beta $. Moreover, ${\mathcal E}_\alpha \otimes 
{\mathcal E}_\beta$ is 
isomorphic to $\mathcal{E}_{\alpha\circ\beta}$. Hence we obtain an homomorphism $\rho_A$ 
from $\mathrm{Out}(A)$ to $\mathrm{Pic}(A)$. 
 Since $A$ is unital, any 
$A$-$A$-equivalence bimodule is a finitely generated projective $A$-module as a right module with a finite basis $\{\xi_i\}_{i=1}^n$. 
Put  $p=(\langle\xi_i,\xi_j\rangle_A)_{ij} \in M_n(A)$. 
Then $p$ is a projection and $\mathcal{E}$ is isomorphic to 
$pA^n$ as right Hilbert $A$-module 
with an isomorphism  of $A$ to  $pM_n(A)p$ as $C^*$-algebra.  

For any unital subring $R$ of $\mathbb R$, we denote by 
$R^{\times}$ (resp. $R^{\times}_{+}$) 
the set of invertible elements (resp. positive invertible elements)
in the ring  $R$. In particular we denote by $\mathbb{R}^{\times}_{+}$ 
the multiplicative group of 
positive real numbers.  

\begin{pro}\label{pro:phom}
Let $A$ be a unital simple $C^*$-algebra with a unique normalized trace $\tau$. Define a map $T: \mathrm{Pic}(A) \rightarrow \mathbb{R}^{\times}_{+}$ by 
$T([\mathcal{E}])=\sum_{i=1}^n\tau (\langle\xi_i,\xi_i\rangle_A)$, 
 where $\{\xi_i\}_{i=1}^n$ is 
a finite basis of $\mathcal{E}$. Then $T([\mathcal{E}])$  
does not depend on the choice of basis and $T$ is well-defined. 
Moreover $T$ is a multiplicative map 
and $T(\mathcal{E}_{id}) = 1$.  
\end{pro}
\begin{proof} Let $\{\xi_i\}_{i=1}^n$ and 
$\{\zeta_k\}_{k=1}^l$ be finite bases of  $\mathcal{E}$. 
Then 

\begin{align*}
\sum_{i=1}^n\tau (\langle\xi_i,\xi_i\rangle_A) 
& =\sum_{i=1}^n\tau (\langle\sum_{k=1}^l\zeta_k 
\langle\zeta_k,\xi_i\rangle_A,\xi_i\rangle_A)=\sum_{i,k=1}^{n,l}\tau (\langle\xi_i,\zeta_k\rangle_A
\langle\zeta_k,\xi_i\rangle_A) \\
& =\sum_{i,k=1}^{n,l}\tau (\langle\zeta_k,
\xi_i\rangle_A\langle\xi_i,\zeta_k\rangle_A) =\sum_{k=1}^l\tau (\langle\zeta_k,\zeta_k\rangle_A).
\end{align*} 
Therefore $T([\mathcal{E}])$  is independent on the choice of basis. 

Let $\mathcal{E}_1$ and $\mathcal{E}_2$ be $A$-$A$-equivalence bimodules with 
bases 
$\{\xi_i\}_{i=1}^n$ and $\{\eta_j\}_{j=1}^m$. 
Suppose that there exists an isomorphism $\Phi$ of $\mathcal{E}_1$ onto $\mathcal{E}_2$. Then  $\{\Phi(\xi_i)\}_{i=1}^n$ is also 
a basis of $\mathcal{E}_2$. Then 
$$
\sum_{i=1}^n\tau (\langle\xi_i,\xi_i\rangle_A)
= \sum_{i=1}^n\tau ( \langle\Phi (\xi_i),\Phi (\xi_i)\rangle_A) 
= \sum_{i=1}^n\tau (\langle\eta_i,\eta_i\rangle_A).
$$
Therefore $T$ is well-defined.  

We shall show that $T$ is multiplicative. 
Let $\mathcal{E}_1$ and $\mathcal{E}_2$ be $A$-$A$-equivalence bimodules with 
bases 
$\{\xi_i\}_{i=1}^n$ and $\{\eta_j\}_{j=1}^m$. 
Then $\{\xi_i\otimes\eta_j\}_{i,j=1}^{n,m}$ is a basis of $\mathcal{E}_1\otimes\mathcal{E}_2$ and 
$$
T([\mathcal{E}_1]\otimes [\mathcal{E}_2])
= \sum_{i,j} \tau (\langle\xi_i\otimes \eta_j,\xi_i\otimes \eta_j \rangle_A) 
=\sum_{j=1}^{m}\tau (\langle\eta_j,
\sum_{i=1}^n\langle\xi_i,\xi_i\rangle_A\eta_j\rangle_A). 
$$ 
Define $\tau^{\prime}(a)=\sum_{j=1}^m
\tau (\langle\eta_j,a\eta_j\rangle_A)$ for $a\in A$. Then we have 
\begin{align*}
\tau^{\prime}(ab) & =
\sum_{j=1}^m\tau (\langle\eta_j,a\sum_{i=1}^m\eta_i\langle\eta_i,b\eta_j\rangle_A)=\sum_{i,j=1}^m
\tau (\langle\eta_j,a\eta_i\rangle_A\langle\eta_i,b\eta_j\rangle_A) \\
& =\sum_{i,j=1}^m\tau (\langle
\eta_i,b\eta_j\rangle_A\langle\eta_j,a\eta_i\rangle_A)=\tau^{\prime}(ba).
\end{align*} 
Therefore $\tau^{\prime}$ is a trace on $A$ and 
$\tau^{\prime}(1)=\sum_{j=1}^m\tau (\langle\eta_j,\eta_j\rangle_A)$. 
Since $\tau$ is the unique normalized trace on $A$, 
$$
\tau^{\prime}(a) = 
\sum_{j=1}^m
\tau (\langle\eta_j,a\eta_j\rangle_A) 
= (\sum_{j=1}^m\tau (\langle\eta_j,\eta_j\rangle_A) \tau (a).
$$ 
Then 
\begin{align*}
& T([\mathcal{E}_1 \otimes_A 
\mathcal{E}_2])
 =\sum_{j=1}^{m}\tau (\langle\eta_j,
\sum_{i=1}^n\langle\xi_i,\xi_i\rangle_A\eta_j\rangle_A)
=\tau^{\prime}(\sum_{i=1}^n\langle\xi_i,\xi_i\rangle_A) \\
& = (\sum_{j=1}^m\tau (\langle\eta_j,\eta_j\rangle_A)
(\tau(\sum_{i=1}^n \langle\xi_i,\xi_i\rangle_A))
=T([\mathcal{E}_1])T([\mathcal{E}_2]). 
\end{align*} 
\end{proof}

\section{Fundamental groups}\label{sec:fun}

Let $A$ be a simple unital $C^*$-algebra with a unique normalized trace $\tau$. We denote by $Tr$ the usual unnormalized trace on $M_n(\mathbb{C})$. 
 Put 
$$
\mathcal{F}(A) :=\{ \tau\otimes Tr(p) \in \mathbb{R}^{\times}_{+}\ | \ 
 p \text{ is a projection in } M_n(A) \text{ such that } pM_n(A)p  \cong A \}. 
$$

\begin{thm}\label{thm:gro}
Let $A$ be a unital simple $C^*$-algebra with a unique normalized trace 
$\tau$. Then $\mathcal{F}(A)$ 
is a multiplicative subgroup of $\mathbb{R}^{\times}_{+}$. Furthermore 
if $A$ is separable, then $\mathcal{F}(A)$ is countable. 
\end{thm}
\begin{proof}
Let $T :  \mathrm{Pic}(A) \rightarrow \mathbb{R}^{\times}_{+}$ be a 
multiplicative map defined in Proposition \ref{pro:phom}. It is 
enough to show that $\mathcal{F}(A) = Im \ T$. 
Let $\mathcal{E}$ be  
an $A$-$A$-equivalence bimodule with 
a finite basis $\{\xi_i\}_{i=1}^n$. 
Define a projection 
 $p=(p_{ij})_{ij} = (\langle\xi_i,\xi_j\rangle_A)_{ij} \in M_n(A)$. 
Then $\mathcal{E}$ is isomorphic to $pA^n$ as right Hilbert $A$-module 
with an isomorphism  of $A$ to  $pM_n(A)p$. Therefore 
$$
T([\mathcal{E}]) = \sum_{i=1}^n\tau (\langle\xi_i,\xi_i\rangle_A) 
=  \sum_{i=1}^n \tau (p_{ii}) =  \tau\otimes Tr(p) . 
$$
Hence  $Im \ T \subset \mathcal{F}(A)$. Conversely  
let $p=(p_{ij})_{ij}$ be a projection in $M_n(A)$ such that $pM_n(A)p$ is isomorphic to $A$. 
Since $A$ is simple, $\mathcal{E} := pA^n$ is an $A$-$A$-equivalence bimodule. 
Thus it is clear that $\mathcal{F}(A) \subset Im \ T$. 
\end{proof}

\begin{Def}
Let $A$ be a unital (separable) simple $C^*$-algebra with a unique normalized trace $\tau$. We call 
$\mathcal{F}(A)$ the fundamental group of $A$, which is a 
multiplicative (countable) subgroup of $\mathbb{R}^{\times}_{+}$. 
\end{Def}

Let $p=(p_{ij})_{ij}$ and $q =(q_{lk})_{lk}$ be projections in $M_n(A)$ and $M_m(A)$ such that $pM_n(A)p\cong A$ and $qM_m(A)q \cong A$. 
Put $s = \tau\otimes Tr(p), t = \tau\otimes Tr(q) \in \mathcal{F}(A)$. 
Then  $st \in \mathcal{F}(A)$ is given by the following projection 
$r \in M_{nm}(A)$: Fix an isomorphism 
$\phi$ of $A$ onto $q M_m(A)q$ given by $\phi(a) = (\phi_{lk}(a))_{lk}$. Put
$r=(\phi_{lk}(p_{ij}))_{(i,l),(j,k)} \in M_{nm}(A)$. Then $st = \tau\otimes Tr(r)$. 

\begin{rem}
Let $M$ be a factor of type $II_1$ and $p$ and $q$ be projections in $M_n(M)$  
such that  $\tau\otimes Tr(p)=\tau\otimes Tr(q)$.  
Then $pM_n(M)p\cong qM_n(M)q$. But in our definition of the 
fundamental group of a unital simple $C^*$-algebra $A$ with a unique trace 
$\tau$, we do {\it not} assume that the $C^*$-algebra $A$ has such a 
property in order to avoid loosing many interesting examples. 
We say that   
the trace $\tau$  on $A$ separates  equivalence classes of projections if 
$\tau\otimes Tr(p)=\tau\otimes Tr(q)$ for projections $p, q \in M_n(A)$ implies that $p$ and $q$ are von Neumann equivalent. 
For example, let $A$ be 
an AF-algebra such that $K_0(A)=\mathbb{Q}\oplus\mathbb{Z}$, 
$K_0(A)_+=\{(q, n)\in K_0(A):q>0 \}\cup \{(0,0) \}$ and $[1]_0=(1,0)$ 
Then  $A$ is a simple unital 
$C^*$-algebra with a unique normalized trace. Take a projection $p\in M_n(A)$ such that 
$[p]_0=(1,1) \in K_0(A)$. Then $\tau\otimes tr(p)=1$ and $pM_n(A)p$ is not isomorphic to 
$A$. The trace $\tau$  on $A$ does not separate equivalence classes of projections.   
\end{rem}

\begin{lem}
Let $A$ and $B$ be  unital simple $C^*$-algebras with unique normalized traces. Then $A\otimes_{\mathrm{min}}B$ is a unital simple $C^*$-algebra with a unique 
normalized trace and  
$$
\mathcal{F}(A)\mathcal{F}(B)\subseteq\mathcal{F}
(A\otimes_{\mathrm{min}}B). 
$$
\end{lem} 
\begin{proof} Obvious.
\end{proof}

We shall show that K-theoretical obstruction enables us to compute 
fundamental groups easily. 
We denote by 
$\tau_*$ the map $K_0(A)\rightarrow \mathbb{R}$ induced by a trace $\tau$ 
on $A$. 

\begin{Def}
Let $E$ be an additive subgroup of $\mathbb{R}$ containing $\mathbb{Z}$. 
Then the {\it inner multiplier group} $IM(E)$ of $E$ is defined by 
$$
IM(E) = \{t \in {\mathbb R}^{\times} \ | t \in E, t^{-1}  \in E, \text{ and } 
        tE = E \}.
$$
Then  $IM(E)$ is a multiplicative subgroup of $\mathbb{R}^{\times}$. 
We call  $IM_+(E) := IM(E) \cap  \mathbb{R}_+$ the 
{\it positive inner multiplier group} of $E$, which is a multiplicative subgroup of $\mathbb{R}^{\times}_+$.   
\end{Def}

If $E$ is a unital subring of  $\mathbb{R}$ , then 
$IM(E)$ (resp. $IM_+(E)$) coincides with the group of 
 invertible elements (resp. positive invertible elements) in the ring $E$.  

\begin{lem}
Let $E$ be an additive subgroup of $\mathbb{R}$ containing $\mathbb{Z}$.  
Consider a multiplicative subgroup  $G = IM_+(E)$ of $\mathbb{R}^{\times}_+$.  
Let 
$$
R_E := \{ t \in \mathbb{R} \ | \ t = \sum_{g \in G} a_gg,  \ a_g \in \mathbb{Z} \text { is zero except finite } g \in G  \} \subset E.
$$
Then $R_E$ is a unital subring of  $\mathbb{R}$ and 
$IM_+(E) = IM_+(R_E) = (R_E)^{\times}_+$. 
\end{lem} 
\begin{proof} 
Since $gE = E$ for any $g \in G$, 
we have $tE = E$  for any $t \in (R_E)^{\times}_+$.  
Thus $IM_+(E) \subset (R_E)^{\times}_+ \subset IM_+(E)$.
\end{proof}

Not all multiplicative subgroups of $\mathbb{R}^{\times}$ arise in this 
manner. For example, let 
$\Gamma = \{ (\sqrt{5})^n \in \mathbb{R}^{\times}_+ \ | \ n \in \mathbb{Z} \}$. Suppose that there were  a unital subring of $R$ of $\mathbb{R}$ such that 
$\Gamma = R^{\times}_{+}$. Since  $R$ contains 
${\mathbb Z} + {\mathbb Z}{\sqrt{5}}$, 
$(\sqrt{5} - 2) \in R^{\times}_{+} = \Gamma$. This is a contradiction. 
By considering $R = R_E$, we know that  $\Gamma$ does not arise as $IM_+(E)$ 
for any additive subgroup of $\mathbb{R}$ containing $\mathbb{Z}$.

\begin{pro}\label{prop:inner-multiplier}
Let $A$ be a simple unital $C^*$-algebra with a unique normalized trace 
$\tau$. Then 
$$
\mathcal{F}(A) \subset IM_+(\tau_*(K_0(A))) \subset 
\tau_*(K_0(A)) \cap \tau_*(K_0(A))^{-1} \cap \mathbb{R}_+.
$$
In particular, $\tau_*(K_0(A))$ is a $\mathbb{Z}[\mathcal{F}(A)]$-module. 
\end{pro} 
\begin{proof}We say that right Hilbert $A$-modules $\mathcal{X}$ and 
$\mathcal{Y}$ with finite bases isomorphic if there exists an right $A$-module isomorphism 
preserving inner product. 
Let $\mathcal{H}(A)$ be the 
set of isomorphic classes $[\mathcal{X}]$ of 
right Hilbert $A$-modules $\mathcal{X}$ with finite basis. 
Define a map $\hat{T} : \mathcal{H}(A) \rightarrow 
\mathbb{R}_{+}$ by 
$\hat{T}([\mathcal{X}])=\sum_{i=1}^n\tau (\langle\xi_i,\xi_i\rangle_A)$, 
 where $\{\xi_i\}_{i=1}^n$ is 
a finite basis of $\mathcal{X}$. Then, as in the proof 
Proposition \ref{pro:phom},  
$\hat{T}([\mathcal{X}])$  
does not depend on the choice of basis and $\hat{T}$ is well-defined. 
Moreover for any  $A$-$A$-equivalence bimodule $\mathcal{E}$, 
$$
\hat{T}([\mathcal{X} \otimes_A \mathcal{E}])
 =T([\mathcal{X}]) \hat{T}([\mathcal{E}]). 
$$
This implies that $\mathcal{F}(A) \subset IM_+(\tau_*(K_0(A))$.
\end{proof} 

\begin{cor}
Let $A$ be a simple unital $C^*$-algebra with a unique normalized trace. Assume that $K_0(A)$ is finitely generated and its torsion free part is 
isomorphic to $\mathbb{Z}$. Then $\mathcal{F}(A)=\{1 \}$. 
\end{cor}
\begin{proof}The assumption implies that $\tau_*(K_0(A))$ is a 
singly generated additive subgroup of $\mathbb{R}$ containing  $\mathbb{Z}$, 
say $\mathbb{Z}\alpha$, 
for some positive real number $\alpha$. 
Since  $\mathcal{F}(A) \cap (0,1)$ is included to a finite set 
$\tau_*(K_0(A)) \cap (0,1) = \mathbb{Z}\alpha \cap (0,1)$, $\mathcal{F}(A)$ is also a finite set. 
This implies that $\mathcal{F}(A)=\{1 \}$.
\end{proof}

We shall show some examples.
\begin{ex}
Let $\mathbb{F}_n$ be a non-abelian free group with $n\geq 2$ generators. Then $C_r^*(\mathbb{F}_n)$ 
is a unital simple $C^*$-algebras with a unique normalized trace. Since $K_0(C_r^*(\mathbb{F}_n))
\cong \mathbb{Z}$, $\mathcal{F}(C_r^*(\mathbb{F}_n))=\{1\}$. 
\end{ex}

\begin{ex}
Let $A=C_r^*(\mathbb{Z}/n\mathbb{Z}*\mathbb{Z}/m\mathbb{Z})$ where $n,m\in\mathbb{N}$ 
such that $(n-1)(m-1)\geq 2$. Then $A$ is a unital simple $C^*$-algebras 
with a unique 
normalized trace. Since $\tau_*(K_0(A))=\frac{1}{\mu (m.n)}\mathbb{Z}$,  
where $\mu (m,n)$ is the lowest common multiple of $m$ and $n$ \cite{cun}, 
$\mathcal{F}(A)=\{1\}$. 
\end{ex}

\begin{ex}
Let $p$ be  a prime number. Consider a UHF algebra  $A=M_{p^\infty}$. 
Then $\tau_*(K_0(A)) = \{\frac{m}{p^n} \ | \ 
m \in \mathbb{Z}, n \in \mathbb{N} \} $ and 
$$
\{p^n:n\in \mathbb{Z}\} \subset \mathcal{F}(A) \subset 
\tau_*(K_0(A)) \cap \tau_*(K_0(A))^{-1} = \{p^n:n\in \mathbb{Z}\}. 
$$
Hence $\mathcal{F}(A)=\{p^n:n\in \mathbb{Z}\}$. 
\end{ex}
\begin{ex}
Let $p$ be  a prime number and $n\geq 2$ an integer. Consider
$B=M_{p^\infty}\otimes C_r^*(\mathbb{F}_n)$. Then $\tau_*(K_0(B)) =  
\tau_*(K_0(M_{p^\infty}\otimes C_r^*(\mathbb{F}_n)))=\tau_*(K_0(M_{p^\infty}))$ 
by \cite{P}(Theorem 3). Therefore 
$$
\{p^n:n\in \mathbb{Z}\} \subset \mathcal{F}(B) \subset 
\tau_*(K_0(B)) \cap \tau_*(K_0(B))^{-1} \cap \mathbb{R}_+ 
= \{p^n:n\in \mathbb{Z}\}. 
$$
Hence $\mathcal{F}(B)=\{p^n:n\in \mathbb{Z}\}$ 
\end{ex}
In a similar way, we have the following proposition. 

\begin{pro}
Let $S$ be a subset of the prime numbers and $G$ the multiplicative 
subgroup of $\mathbb{R}^{\times}_+$ generated by $S$. Then 
there exist a separable  nuclear simple $C^*$-algebra $A$ with a 
unique trace and a  separable  non-nuclear simple $C^*$-algebra $B$ 
with a unique trace such that $\mathcal{F}(A) = \mathcal{F}(B) = G$. 
\end{pro}
\begin{proof}
Let $m =\prod _{p\in S} p^{\infty}$ be a supernatural number and 
$A$ the corresponding UHF algebra. Then   
$$
G \subset \mathcal{F}(A) \subset 
\tau_*(K_0(A)) \cap \tau_*(K_0(A))^{-1} = G. 
$$
Hence $\mathcal{F}(A)= G$. Put $B=A \otimes C_r^*(\mathbb{F}_2)$. 
Since  $\tau_*(K_0(B)) =  \tau_*(K_0(A)))$, 
$$
G \subset \mathcal{F}(B) \subset 
\tau_*(K_0(B)) \cap \tau_*(K_0(B))^{-1} \cap \mathbb{R}_+ 
= G. 
$$
Hence $\mathcal{F}(A)= G$. 
\end{proof}

The following theorem shows that  many countable subgroups of 
$\mathbb{R}^{\times}_{+}$ 
can be fundamental 
groups.  

\begin{thm} \label{thm:subring}
Let $R$ be a unital countable subring of $\mathbb{R}$. 
Then there exist a 
separable,  unital simple, nuclear  $C^*$-algebra $A$ and a
non-nuclear $C^*$-algebra $B$ with unique trace such that 
$\mathcal{F}(A)=\mathcal{F}(B)=R^\times_+$. 
\end{thm}
\begin{proof}
Since$(R,R_+,1)$ is a simple unperforated ordered group that has the 
Riesz interpolation property, there exists a simple $AF$-algebra $A$ such that 
$(K_0(A),K_0(A)_+,[1_A]_0)=(R,R_+,1)$ by \cite{EHS}. Then $A$ has a unique 
normalized trace $\tau$ and $\tau_*(K_0(A))=R$. 
By Proposition \ref{prop:inner-multiplier}, 
$\mathcal{F}(A)\subset R^\times_+$. Conversely let  $\lambda\in R^\times_+$. 
Then there exists a projection $p \in M_n(A)$ for some $n\in\mathbb{N}$ such 
that $\tau\otimes Tr(p)=\lambda$. Define a additive 
homomorphism $\varphi :R \rightarrow R$ by $\varphi (r) =\lambda r$ for 
$r \in R$ Then $\varphi$ is an order 
isomorphism with $\varphi (1)=\lambda$. Since 
$(K_0(pM_n(A)p),K_0(pM_n(A)p)_+,[p]_0) =
(R,R_+,\lambda)$, there exists an isomorphism $f:pM_n(A)p \rightarrow A$ 
with $f_* = \varphi$ by \cite{E}. Therefore $\lambda\in\mathcal{F}(A)$. 

Put $B=A\otimes C_r(\mathbb{F}_2)$. Then $B$ is a simple separable non-nuclear $C^*$-algebra with unique normalized trace. 
Since  $\tau_*(K_0(B)) =  \tau_*(K_0(A))) = R$, 
$$
R^\times_+ \subset \mathcal{F}(B) \subset IM_+(\tau_*(K_0(B))) = 
IM_+(R) =R^\times_+ . 
$$
Hence $\mathcal{F}(A)= R^\times_+$. 
\end{proof}

\begin{rem}
By the proof above, if $A$ a simple unital separable AF-algebra with a unique normalized trace 
which separates equivalence classes of projections, then 
$\mathcal{F}(A\otimes C_r^*(\mathbb{F}_n))=\mathcal{F}(A)$. 

\end{rem}

\begin{cor}
For any countable subset $S \subseteq\mathbb{R}^\times_+$, there exist a separable, 
simple,  nuclear $C^*$-algebra $A$ and non-nuclear $C^*$-algebra $B$ with 
unique trace such that $S\subseteq\mathcal{F}(A)$ and 
$S\subseteq\mathcal{F}(B)$. 
\end{cor}
\begin{proof}
Let $G$ be the multiplicative subgroup of $\mathbb{R}^\times_+$ generated 
by $S$. Let
$$
R= \{ t \in \mathbb{R} \ | \ t = \sum_{g \in G} a_gg,  \ a_g \in \mathbb{Z} \text { is zero except finite } g \in G  \} 
$$
Then $R$ is a unital subring of $\mathbb{R}$ with 
$S \subset G \subset  R^\times_+$.
Applying the above Theorem \ref{thm:subring} for $R$ completes the proof. 
\end{proof}

If a $C^*$-algebra $A$ belongs to the class covered by Elliott's 
classification theorem, then the fundamental group ${\mathcal F}(A)$ 
can be computed rather easily.

\begin{pro}\label{pro:AT}
Let $A$ be a unital simple separable A$\mathbb{T}$-algebra of real rank zero with a unique normalized trace. 
Assume that $\tau_* : K_0(A) \rightarrow \tau_*(K_0(A))$ is 
an order isomorphism. 
Then  $\mathcal{F}(A) = IM_+(\tau_*(K_0(A)))$.  
\end{pro}
\begin{proof} Let $E = \tau_*(K_0(A))$.  By Proposition \ref{prop:inner-multiplier}, 
we have 
$\mathcal{F}(A)\subseteq IM_+(E)$. Conversely 
$\lambda\in IM_+(E)$, then there exist a projection $p$ 
in $M_n(A)$ for some $n\in\mathbb{N}$ such that $\tau\otimes tr(p)=\lambda$. Define an additive  
homomorphism $\varphi :E \rightarrow E$ by $\varphi (r) =\lambda r$. 
Then $\varphi$ is an order isomorphism and $\varphi (1)=\lambda$. 
It is clear that  $(K_0(pM_n(A)p),K_0(pM_n(A)p)_+,[p]_0) = (E,E_+,\lambda)$ and  $K_1(pM_n(A)p)$ is isomorphic 
to $K_1(A)$. 
By a classification 
theorem of $A\mathbb{T}$-algebras of real rank zero by Elliott 
\cite{Ell}, there exists an isomorphism of $A$ to $pM_n(A)p$. 
Therefore $\lambda\in\mathcal{F}(A)$. 
\end{proof}

We shall compute the fundamental group of irrational rotation algebras 
$A_\theta$ 
for irrational numbers $\theta$. 
We know that  $\tau_*(K_0(A_\theta ))=\mathbb{Z}+\mathbb{Z}\theta$ 
by \cite{PV},\cite{Rie}. 
But we should note that $\mathbb{Z}+\mathbb{Z}\theta$ 
is not a ring in general. 
 
\begin{cor}\label{cor:irrational}
Let $\theta$ be an irrational number. 
Then $\mathcal{F}(A_\theta ) = 
IM_+(\mathbb{Z}+\mathbb{Z}\theta )$.  
More precisely, 
if $\theta$ is not a 
quadratic number, then 
$\mathcal{F}(A_\theta )=\{1\}$.   
If $\theta$ is a quadratic number, 
then $\mathcal{F}(A_\theta ) 
= \{\epsilon _0^n  \in R^\times_+ \ | \ n \in \mathbb{Z} \}$ 
is singly generated and the generator $\epsilon _0$ is given by 
$\epsilon _0 =\frac{t+u\sqrt{D_\theta}}{2}$, 
where 
$D_\theta$ is the  discriminant of $\theta$ and $t,u$ are the positive 
integers satisfying one of the Pell 
equations $t^2-D_\theta u^2=\pm 4$ and $\frac{t+u\sqrt{D_\theta}}{2}$ exceeds $1$ 
and is least. 
\end{cor}
\begin{proof}
By Proportion \ref{pro:AT}, we have 
$
\mathcal{F}(A) = IM_+(\tau_*(K_0(A))) 
=IM_+(\mathbb{Z}+\mathbb{Z}\theta). 
$
If $\theta$ is not a quadratic number, 
then 
$$
\mathcal{F}(A_\theta) \subset 
(\mathbb{Z}+\mathbb{Z}\theta) \cap (\mathbb{Z}+\mathbb{Z}\theta)^{-1} 
= \{1\}.  
$$
Assume that $\theta$ is a quadratic 
algebraic integer. 
Then $\mathbb{Z}+\mathbb{Z}\theta$ is a unital ring and 
$\mathcal{F}(A_\theta ) = 
IM_+(\mathbb{Z}+\mathbb{Z}\theta ) = (\mathbb{Z}+\mathbb{Z}\theta )^{\times}_+$. 
Let $\epsilon\in (\mathbb{Z}+\mathbb{Z}\theta )^{\times}_+$. 
Then there exist integers $x,y$ such that $\epsilon =\frac{x+y\sqrt{D_\theta}}{2}$ and 
satisfying one of the Pell equations $x^2-D_\theta y^2=\pm 4$. 
If $\frac{t+u\sqrt{D_\theta}}{2}>1$, then $t>0$ and $u>0$. 
Therefore 
$(\mathbb{Z}+\mathbb{Z}\theta )^{\times}_+=
 \{\epsilon _0^n  \in  R^\times_+ \ | \ n \in \mathbb{Z} \}$ 
is singly generated and the generator $\epsilon _0$ is given by 
$\epsilon _0 = \frac{t+u\sqrt{D_\theta}}{2}$ as above. 
(We refer the reader to \cite{JW} or \cite{BV} for details.)

Let $\theta$ be a quadratic number and  not an 
algebraic integer with $k\theta^2+l\theta +m=0$ for some integer
$k>1, l,m$ with $gcd(k,l,m)=1$. Then 
$\mathbb{Z}+\mathbb{Z}\theta$ is not a ring. 
Since  $k\theta$ is a quadratic integer, 
$\mathbb{Z}+\mathbb{Z}k\theta$ is a unital subring. It is easy to 
verify that 
$(\mathbb{Z}+\mathbb{Z}k\theta )^\times 
= IM(\mathbb{Z}+\mathbb{Z}\theta)$.    
Since $D_\theta =D_{k\theta}$, the computation of $\mathcal{F}(A_\theta )$
is reduced to the case of algebraic integer.

\end{proof}

\begin{ex}
Let $\theta =\sqrt{3}$. Then 
$\tau_*(A_\theta) = \mathbb{Z}+\mathbb{Z}\sqrt{3}$ is the  ring of integers of 
a quadratic field $\mathbb{Q}(\sqrt{3})$ and $\mathcal{F}(A)=\{(2+\sqrt{3})^n:n\in\mathbb{
Z}\}$. 
\end{ex}
\begin{ex}
Let $\theta =\frac{-1+\sqrt{5}}{2}$. Then 
$\tau_*(A_\theta) = \mathbb{Z}+\mathbb{Z}\frac{-1+\sqrt{5}}{2}$ 
is the  ring 
of integers of a quadratic field 
$\mathbb{Q}(\sqrt{5})$ and $\mathcal{F}(A)=\{(\frac{1+\sqrt{5}}{2})^n:n 
\in \mathbb{Z}\}$. 
\end{ex}
\begin{ex}
Let $\theta =\sqrt{5}$. Then 
$\tau_*(A_\theta) = \mathbb{Z}+\mathbb{Z}{\sqrt{5}}$ is a ring and 
$\mathcal{F}(A)=\{(\sqrt{5}+2)^n:n\in\mathbb{Z}\}$. 
\end{ex}
\begin{ex}
Let $\theta =\frac{1}{\sqrt{5}}$. Then 
$\tau_*(A_\theta) = \mathbb{Z}+\mathbb{Z}\frac{1}{\sqrt{5}}$ is not a ring 
and 
$\mathcal{F}(A)=\{(\sqrt{5}+2)^n:n\in\mathbb{Z}\}$. 
\end{ex}

\begin{cor}
Let $\theta$ be an irrational number. 
Then 
$$
\mathcal{F}(A_{\theta}  \otimes C_r^*(\mathbb{F}_n)) = 
 \mathcal{F}(A_\theta) = IM_+(\mathbb{Z}+\mathbb{Z}\theta ).
$$  
\end{cor}
\begin{proof}
Put $B=A_{\theta} \otimes C_r^*(\mathbb{F}_n)$. 
Since  $\tau_*(K_0(B)) =  \tau_*(K_0(A_{\theta})))$ by \cite{P}(Theorem 3), 
$$
 \mathcal{F}(A_\theta)
\subset \mathcal{F}(B) \subset 
IM_+(\tau_*(K_0(B))) = IM_+(\tau_*(K_0(A_\theta))) 
= \mathcal{F}(A_\theta).
$$
\end{proof}

\begin{rem}
Let $A$ be a unital simple separable AH-algebra of slow dimension growth and of real rank zero with the 
unique normalized trace. Assume that $\tau_*(K_0(A))$ is a free abelian group.  Then $\mathcal{F}(A) = IM_+(\tau_*(K_0(A)))$.  
In fact,  
let $E=\tau_*(K_0(A))$. By Proposition \ref{prop:inner-multiplier}, 
$\mathcal{F}(A)\subset IM_+(\tau_*(K_0(A))) $. 
Since $\tau_*(K_0(A))$ is a free abelian group, $K_0(A)$ is isomorphic to $E\oplus ker(\tau_*)$. 
The positive cone $K_0(A)_+$ is the set $\{(r_1,r_2):r_1>0,r_2\in ker(\tau_*)\}\cup\{(0,0)\}$ 
because the normalized trace on $A$ is unique. By classification theorem, 
we have $IM_+(\tau_*(K_0(A))) \subset \mathcal{F}(A)$ by a similar argument 
as in Proposition \ref{pro:AT}. 
\end{rem}

We shall consider a  relation between  fundamental groups and  Picard groups 
of $C^*$-algebras. 

\begin{pro}
Let $A$ be a unital simple $C^*$-algebra with a unique normalized trace $\tau$. Assume that 
the normalized trace on $A$ separates equivalence classes of projections. Then we have the 
following exact sequence. 
\[\begin{CD}
      {1} @>>> \mathrm{Out}(A) @>\rho_A>> \mathrm{Pic}(A) @>T>> \mathcal{F}(A)
 @>>> {1} \end{CD}. \]
\end{pro}
\begin{proof}
It is clear that $\rho_A$ is one to one, $T$ is onto and 
$Im \ \rho_A \subset Ker(T)$. We shall show that 
$Ker(T) \subset Im \ \rho_A$ 
Let $p$ be a projection in $M_n(A)$ with an isomorphism  
$\varphi : A \rightarrow pM_n(A)p$.  
Consider the corresponding  equivalence bimodule 
$\mathcal{E} = pM_n(A)p$.  
Assume that $\rho_A([\mathcal{E}]) = \tau\otimes tr(p)=1$. 
Since the normalized trace on $A$ separates equivalence classes of projections,
there exists a partial isometry $w$ such that $p =w^*w$ and 
$1\otimes e_{11} = ww^*$ in $M_n(A)$ 
where $e_{11}$ is a rank one projection in $M_n(\mathbb{C})$. 
Then there exists an automorphism $\alpha$ of $A$ such that 
$w\varphi(a)w^* = \alpha(a) \otimes e_{11}$ 
Hence $[\mathcal{E}] = [\mathcal{E}_{\alpha}] \in Im \rho_A$.  
\end{proof}
Let $\theta$ be a quadratic number with $k\theta^2+l\theta +m=0$. 
Then $\mathrm{Pic}(A_\theta )$ is isomorphic to $\mathrm{Out}(A_\theta)\rtimes\mathbb{Z}$ 
by Corollary \ref{cor:irrational}. 
In \cite{kod2}, Kodaka showed this fact by considering automorphisms of 
$A_\theta\otimes K(\mathcal{H})$ and showing that 
$\{c\theta +d;\frac{a\theta +b}{c\theta +d}=\theta ,a,b,c,d\in\Bbb{Z},ad-bc=\pm 1 \}$ 
is singly generated (see also \cite{N}). 
Our argument is different from Kodaka's argument. But there exists a connection between 
them. In fact, it is known that 
$$\{c\theta +d;\frac{a\theta +b}{c\theta +d}=\theta ,a,b,c,d\in\Bbb{Z},ad-bc=\pm 1 \}
=(\mathbb{Z}+\mathbb{Z}k\theta )^\times .
$$

Recall that 
the fundamental group of a $II_1$-factor $M$ is equal to the set of trace-scaling constants for 
automorphisms of $M\otimes B(\mathcal{H})$. We have a similar fact as discussed by Kodaka in \cite{kod3}. 
Let $K(\mathcal{H})$ be the $C^*$-algebra of all compact operators on a countably infinite dimensional Hilbert space $\mathcal{H}$.  For $x\in (A\otimes K(\mathcal{H}))_+$, set 
$\hat{\tau} (x)=\sup\{\tau \otimes Tr(y):y\in \cup _n M_n(A), y\leq x\}$. Define $\mathcal{M}_\tau^+=
\{x\geq 0: \hat{\tau}  (x)<\infty \}$ and $\mathcal{M}_\tau =\mathrm{span} \mathcal{M}_\tau^+$. 
Then $\hat{\tau}$ is a densely defined (with the domain $\mathcal{M}_\tau $) 
lower semi-continuous trace 
on $A\otimes K(\mathcal{H})$. Since the normalize trace on $A$ is unique, the lower semi-continuous 
densely defined trace on $A\otimes K(\mathcal{H})$ is unique up to constant multiple. 
We define the set of trace-scaling constants for automorphisms:  
$$
\mathfrak{S}(A)
:= \{ \lambda \in \mathbb{R}^{\times}_+ \ | \ 
\hat{\tau} \circ \alpha (x)=\lambda \hat{\tau}(x) \text{ for } 
 x \in \mathcal{M}_\tau \}. 
$$
It is clear  that $\mathfrak{S} (A)$ is a multiplicative 
subgroup of $\mathbb{R}^{\times}_+$. 
\begin{pro}
Let $A$ be a unital simple $C^*$-algebra with a unique normalized 
trace $\tau$. 
Then $\mathcal{F}(A)=\mathfrak{S}(A)$. 
\end{pro}
\begin{proof}
Let $\lambda \in \mathfrak{S}(A)$, then there exists an automorphism of $A\otimes K(\mathcal{H})$ 
such that $\hat{\tau} \circ \alpha(x)=\lambda \hat{\tau}(x)$ for $x\in\mathcal{M}_\tau$. Let $e$ be a rank one 
projection in $K(\mathcal{H})$. Then 
$I \otimes e \in \mathcal{M}_\tau$ and 
$(I\otimes e)(A\otimes K(\mathcal{H}))(I \otimes e)$ is 
isomorphic to $A$. We have that 
$\alpha (I \otimes e)(A\otimes K(\mathcal{H}))\alpha (I \otimes e)$ is 
isomorphic to $A$ and $\hat{\tau}(\alpha (I \otimes e))=\lambda$.
Therefore $\lambda\in\mathcal{F}(A)$.  

Conversely, let $\lambda\in\mathcal{F}(A)$. 
There exists an projection $p \in M_n(A)$ with an isomorphism 
$\varphi : A \rightarrow pM_n(A)p$  such that 
$\tau \otimes Tr (p) = \lambda$. Since $p$ is a full projection, 
there exists a partial isometry 
$w \in M(M_n(A) \otimes K(\mathcal{H}))$ 
such that $w^*w = I\otimes I$ and $ww^* = p \otimes I$ by 
Brown \cite{B}.  Let 
$\psi : A \otimes M_n(\mathbb{C}) \otimes K(\mathcal{H}) 
\rightarrow A \otimes K(\mathcal{H})$ be an isomorphism which 
induces the  identity on the $K_0$-group. 
Define $ \alpha := \psi \circ (ad \ w^*) 
\circ (\varphi \otimes I) \in Aut(A \otimes K(\mathcal{H}))$
Then $\hat{\tau} (\alpha(I \otimes e_{11})) = 
\hat{\tau} (p \otimes e_{11}) = \lambda$. Therefore 
$\lambda\in  \mathfrak{S}(A)$.     

\end{proof}

Finally we state a direct relation between the fundamental group 
of $C^*$-algebras and that of von Neumann algebras.   

\begin{pro}
Let $A$ be a unital infinite-dimensional simple $C^*$-algebras with a unique normalized 
trace $\tau$. Consider the GNS representation 
$\pi_{\tau} : A \rightarrow B(H_{\tau})$ and the associated 
factor $\pi_{\tau}(A)''$ of type $II_1$. 
Then $\mathcal{F} (A) \subset \mathcal{F} (\pi_{\tau}(A)'')  $. 
In particular, if $\mathcal{F} (\pi_{\tau}(A)'') = \{ 1 \}$, then 
$\mathcal{F}(A) = \{1 \}$. 
\end{pro}
\begin{proof}
Obvious.
\end{proof}
\begin{ex}
Let $S_{\infty}$ be the group of finite permutations of countable infinite set, and let 
$G=S_{\infty}*(\mathbb{Z}^2\rtimes SL(2,\mathbb{Z}))$. Then $C_r^*(G)$ is a unital simple 
$C^*$-algebras with a unique normalized trace. Since $\pi_\tau(A)''$ is isomorphic to 
$L(S_{\infty})*L(\mathbb{Z}^2\rtimes SL(2,\mathbb{Z}))$, $\mathcal{F}(C_r^*(G))=\{1\}$ by 
\cite{IPP}(Corollary 0.3.) and the proposition above. 
\end{ex}

\end{document}